\documentclass{article}

\usepackage{array}
\usepackage{amsmath}
\usepackage{amssymb}
\usepackage{amsfonts}
\usepackage{theorem}
\newtheorem{theorem}{Theorem}
\usepackage{latexsym}

\newcommand\RR{\mathbb{R}}
\begin{document}

\title{A Simple Rule to find a Basic Feasible Solution}
\author{R. Grothmann, Kath. Univ. Eichst\"att-Ingolstadt}

\maketitle

\setlength{\parskip}{10pt plus5pt minus2pt}
\setlength{\parindent}{0pt}

\begin{abstract}
This short note provides and proves an easy algorithm
to find a basic feasible solution for the Simplex Algorithm. The
method uses a rule similar to Bland's rule for the initial phase of
the algorithm.
\end{abstract}


\section{Introduction}

We want to solve a linear optimization problem in the so-called
standard form.
	\[
		\text{Minimize $c^T x$ provided $Ax=b$ and $x \ge
		0$.}           \tag{P}
	\]
Here, $A$ is a $m \times n$ real matrix, $b \in \RR^m$, $c \in
\RR^n$, $x \in \RR^n$. We may assume $\text{rank }(A)=m$. A basic
solution of the system is an $x \in \RR^n$ connected to a tupel
$J$ of column indices
	\[
		1 \le j_1 < \ldots < j_m \le n
	\]
such that the matrix $A_J$ consisting of the columns in $J$ of $A$
is regular (called a column base), and
	\[
		x_J = (A_J)^{-1} b, \qquad \text{$x_j=0$ for $j \notin
		J$}.
	\]
Then clearly $Ax=b$.

Such a basic solution $x$ can be determined easily with the
Gau\ss-algorithm. If the corresponding solution satisfies $x \ge 0$
then $x$ is feasible and the Simplex algorithm can be started to find
an optimal feasible solution, or to find that the problem is
unbounded.

If $x$ is not feasible the standard method of finding a basic feasible
solution is to extend the problem~(P) to
	\[
		\text{Minimize $\sum_{i=1}^m t_i$ provided
		$Ax+t=b$ and $x \ge 0$, $t \ge 0$.}	\tag{PE}	
	\]
Since equalities can be multiplied by $-1$ we can assume $b \ge 0$
in~(P). So (PE) has the basic feasible solution $x=0$ and $t=b$, and
the corresponding column base $J$ contains the columns of the $t_i$. 

Moreover, (P) has a feasible point if and only if (PE) has the minimal
value~$0$. The columns $J$ of the optimal basic solution of (PE) can be
easily used to find a basic feasible solution (P).

This extension is proposed in many books. Examples are the
introductions by Dantzig and Tapa \cite{dantzig}, Alevras and Padberg
\cite{alevras}, Kosmol \cite{kosmol}, and of course the Wikipedia
articles about the Simplex algorithm.

However, this extension is not at all necessary, and there is a much
easier way to find a basic feasible solution using a rule similar to
Bland's rule (see \cite{dantzig}). The purpose of this paper is to
propose such a rule for the initial phase, and to prove that the
algorithm succeeds after a finite number of steps.


\section{The Algorithm}

Starting with a basic solution $x$ as above, let us assume that $x \ge
0$ does not hold. 

Let us denote the scheme after the Gau\ss-Algorithm as
	\[
		\tilde A = A_J^{-1} A, \quad \tilde b = A_J^{-1} b.
	\]
Note $x_J = \tilde b$. Note that we still assume that $J$ is sorted,
so that the unit vectors in $I_m$ appear in $\tilde A$ in the correct
order. By the way, the target function does not matter for the initial
phase.

We then proceed as follows. 

\begin{enumerate}

\item We determine the minimal $i$, such that $x_{j_i} = \tilde b_i
<0$.

\item In the $i$-th row of $\tilde A$ we determine a minimal $j$ such
that $\tilde a_{i,j}<0$.

\item We then exchange the column $j_i \in J$ for $j$.

\end{enumerate}

Note that the problem cannot have a feasible point if the second step
fails since
	\[
		Ax = b \quad\Longleftrightarrow\quad \tilde A x = A_J^{-1} A x
		= A_J^{-1} b = \tilde b.
	\]
In other words, the Gau\ss-algorithm transforms the system $Ax=b$ to
the equivalent system $\tilde A x = \tilde b$. But $\tilde A x =
\tilde b$ is not possible if the $i$-th row of $\tilde A$ is
non-negative and $b_i<0$, due to the restriction $x \ge 0$.

We continue the three steps until either $x \ge 0$, or we see that
there is no feasible solution.

\begin{theorem}
The Algorithm as described above ends after a finite number of steps,
with either a feasible basic solution, or a failure of the second
step showing that there is no feasible solution.
\end{theorem}

\noindent\textbf{Proof:} 

If the algorithm does not terminate we get a cyclic repetition of
column bases
	\[
		J_1 \to J_2 \to \ldots \to J_k \to J_1.
	\]
In these tupels there is a maximal column index $j_u$ that is at one
point inserted and at another point removed in the cycle. There may be
larger columns
	\[
		j_{u+1} < \ldots < j_m
	\]
which are in all column bases of the cycle.

Now there is one of the column bases $J_a$ in the cycle where the
$j_u$ is going to be inserted in the next step. This can only happen if
	\[
		\tilde a_{i,1},\ldots,\tilde a_{i,j_u-1} \ge 0, \quad
		\tilde a_{i,j_u}<0, \quad \tilde b_i<0
	\]
in some row $i$ for the scheme $\tilde A x = \tilde b$ corresponding
to $J_a$. The columns $j_{u+1},\ldots,j_m$ of each scheme $\tilde A$
in the cycle contain the unit vectors $e_{u+1},\ldots,e_m$, and thus
	\[
		\tilde a_{i,j_u+1}=\ldots=\tilde a_{i,j_m}=0.
	\]
We deduce that there is no $x \in \RR^n$
such that
	\[
		x_j = 0 \qquad\text{for $j \ne J$ and $j > j_u$}
	\]
and
	\[
		x_{j_1},\ldots,x_{j_u-1} \ge 0, \quad x_{j_u}<0
	\]
which solves our system of equations, no matter what values we assign
to the variables $x_{j_{u+1}},\ldots,x_{j_m}$.

Then there is another basis $J_b$ in the cycle where the column $j_u$
is going to be removed in the next step. By our rule in step~1, this
can only happen if
	\[
		\tilde b_1, \ldots, \tilde b_{u-1} \ge 0, \quad \tilde
		b_u < 0.
	\]
where $\tilde b$ is the right hand side in our scheme corresponding
to $J_b$. The corresponding basic solution $x$ now satisfies
	\[
		x_j = 0 \qquad\text{for $j \ne J$}
	\]
and
	\[
		x_{j_1}=\tilde b_1,\ldots,x_{j_u-1}=\tilde b_{u-1} \ge 0,
		\quad x_{j_u}=\tilde b_{u}<0.
	\]

Since the schemes are equivalent this is a contradiction. We conclude
that a cycle cannot happen, and thus the algorithm must end.

\hfill{q.e.d.}

It can be shown that we can remove all columns and rows that are never
used for a Pivot element $\tilde a_{i,j}$ during the scheme, so
effectively $j_u=m$. However, this additional step requires some
additional arguments, and we found the proof above easier to
understand.

If the algorithm is computed with the Gau\ss-Algorithmus without
sorting the base columns, we have to follow the following rule: Search
for a Pivot $(i,j)$ such that
	\[
		\tilde a_{i,j}< 0, \quad \tilde b_i<0
	\]
in the scheme, and such that it removes a column with minimal column
from the base, and adds a column with minimal index. Since the
selection of $i$ determines the column that will be removed this must
be done first. Then $j$ can simply be chosen as the minimal index that
satisfies $a_{i,j}<0$.

\end{document}